\documentclass[reqno]{amsart}
\usepackage{amssymb, amsmath, amsthm, amscd, mathrsfs}

\usepackage{paralist}
\usepackage{cite}
\usepackage{bigints}
\usepackage{dsfont}
\usepackage{empheq}
\usepackage{amssymb}
\usepackage{cases}
\usepackage{enumitem}
\allowdisplaybreaks
\usepackage{bbm}
\usepackage{caption}
\usepackage{booktabs}
\usepackage{float}
\usepackage[ruled,vlined,linesnumbered,noresetcount]{algorithm2e}

\usepackage{hyperref}
\numberwithin{equation}{section}

\theoremstyle{plain}
\newtheorem{theorem}{Theorem}

\newtheorem{lemma}{Lemma}
\newtheorem{proposition}{Proposition}
\theoremstyle{definition}

\newtheorem{remark}{Remark}

\newcommand{\vertiii}[1]{{\left\vert\kern-0.25ex\left\vert\kern-0.25ex\left\vert #1 
    \right\vert\kern-0.25ex\right\vert\kern-0.25ex\right\vert}}

\def \d {\mathrm{d}}

\title[Logarithmic convexity of evolution equations] 
      {Logarithmic convexity of evolution equations and application to inverse problems}      
 
\author{Salah-Eddine Chorfi}

\address{S. E. Chorfi, Cadi Ayyad University, UCA, Faculty of Sciences Semlalia, Laboratory of Mathematics, Modeling and Automatic Systems, B.P. 2390, Marrakesh, Morocco}
\email{s.chorfi@uca.ac.ma}

\makeatletter 
\@namedef{subjclassname@2020}{%
  \textup{2020} Mathematics Subject Classification}
\makeatother

 
\begin{document}
\begin{abstract}
We review some results on the logarithmic convexity for evolution equations, a well-known method in inverse and ill-posed problems. We start with the classical case of self-adjoint operators. Then, we analyze the case of analytic semigroups. In this general case, we give an explicit estimate, which will be used to study inverse problems for initial data recovery. We illustrate our abstract result by an application to the Ornstein-Uhlenbeck equations. We discuss both analytic and non-analytic semigroups. We conclude with recent results for the time-fractional evolution equations with the Caputo derivative of order $0<\alpha <1$. We start with symmetric evolution equations. Then, we show that the results extend to the non-symmetric case for diffusion equations, provided that a gradient vector field generates the drift coefficient. Finally, some open problems will be mentioned.
\end{abstract}
\maketitle

\section{Introduction and motivation}
\label{sec:1}

Let $T>0$ be a final fixed time and $A: D(A) \subset H \rightarrow H$ the generator of a $C_0$-semigroup $\left(\mathrm{e}^{t A}\right)_{t \geq 0}$ on a Hilbert space $(H,\langle \cdot, \cdot\rangle)$. We study the inverse problem of determining a class of initial data $u_0$ in an abstract problem
\begin{empheq}[left = \empheqlbrace]{alignat=2}\label{absy}
\begin{aligned}
& u'(t)=A u(t), \quad t \in (0,T], \\
& u(0)=u_0 \in H, 
\end{aligned}
\end{empheq}
from the observation
$$v(t)=\mathbf{C}u(t), \qquad t\in (0,T],$$
where $\mathbf{C}\in \mathcal{L}(V,Y)$ is an observation operator for $(\mathrm{e}^{t A})_{t \ge 0}$ ($V$ and $Y$ are given Hilbert spaces).
\bigskip

\noindent\textbf{Example.} A typical example is given by the heat equation in a bounded domain $\Omega$ of the Euclidean space:
\begin{empheq}[left = \empheqlbrace]{alignat=2}
\begin{aligned}
& \partial_t u(t,x)=\Delta u(t,x), &&\quad (t,x) \in (0,T)\times\Omega, \\
& u=0, &&\quad (t,x) \in (0,T)\times\partial \Omega,\\
& u(0,\cdot)=u_0 \in L^2(\Omega), &&\quad \text{in } \Omega,\nonumber
\end{aligned}
\end{empheq}
with either distributed (interior) or boundary observations:
$$\mathbf{C}u=\mathds{1}_\omega u\qquad \text{ or } \qquad \mathbf{C}u=\partial_n u\rvert_{\Gamma_0},$$
where $\omega \subset \Omega$ and $\Gamma_0\subset \partial \Omega$ are nonempty and relatively open subsets, $\mathds{1}_\omega$ is the indicator function of $\omega$, and $\partial_n$ is the normal derivative with the outer unit normal vector field $n$.

\subsection{Well-posedness of the inverse problem}
We are concerned with the well-posedness of the inverse problem:
\begin{itemize}
    \item[(i)] Uniqueness: Does $\mathbf{C}u=0$ imply $u_0=0$?
    \item[(ii)] Stability: Is it possible to estimate $\|u_0\|$ by a suitable norm of $\mathbf{C}u$?
    \item[(iii)] Numerical reconstruction: Can we design efficient algorithms to reconstruct $u_0$?
\end{itemize}
First, we aim to show a logarithmic stability estimate for a class of initial data:
$$\|u_0\| \le \frac{C}{\left| \log \|\mathbf{C}u\|_{L^2(0,T;Y)}\right|^\alpha},$$
for some $\alpha \in (0,1]$. The general idea can be summarized as follows:
\begin{center}
Observability inequality  \;+\;  Logarithmic convexity \; $\Longrightarrow$ \; Logarithmic stability.
\end{center}

Let us start by recalling the notion of observability. The system \eqref{absy} is final state observable in time $T$ if there exists a constant $\kappa_{T}>0$:
\begin{equation}\label{obs}
\forall u_{0} \in D(A), \qquad\left\|u(T)\right\|^{2} \leq \kappa_{T} \int_{0}^{T}\left\|\mathrm{C} u(t)\right\|_{Y}^{2} \mathrm{d} t.
\end{equation}

Next, we focus on the logarithmic convexity.
\section{Self-adjoint case}
\begin{lemma}[Agmon-Nirenberg (1963)]
Assume that $A$ is self-adjoint. The solution $u$ of \eqref{absy} satisfies
\begin{equation*}
\|u(t)\| \leq \|u_0\|^{1- \frac{t}{T}} \|u(T)\|^{\frac{t}{T}}
\end{equation*}
for all $0\leq t\leq T$.
\end{lemma}
The key ideas of the proof rely on differentiating $\log\|u(t)\|$ twice with respect to $t$, using symmetry and the Cauchy-Schwarz inequality. See the proof of \cite[Theorem 1.3]{GT11}.

\begin{remark} A few remarks are in order:
\begin{itemize}
\item A logarithmic convexity estimate ($K>0$ is constant): $$\|u(t)\| \leq K \|u_0\|^{1- \frac{t}{T}} \|u(T)\|^{\frac{t}{T}} \quad\text{for all } 0\leq t\leq T$$ implies the backward uniqueness property of the solution, i.e., if $u(T)=0$, then $u_0=0$.
\item A well-posed problem need not satisfy logarithmic convexity. A counter-example is given by $u_t + u_x=0, \; (t,x) \in (0,T)\times (0,1),$ $u(t,0)=0, \; (t,x) \in (0,T)\times (0,1)$, $u(0,\cdot)=u_0 \in L^2(0,1)$.
\end{itemize}
We refer to the book \cite{Pa75} and the references therein for more details on the logarithmic convexity method.
\end{remark}

\subsection*{Stability estimate in the self-adjoint case}
Assume that $A$ is self-adjoint and dissipative. For fixed $\epsilon \in(0,1)$ and $M>0$, consider the admissible set of initial data
$$
\mathcal{I}_{\epsilon, M}:=\left\{u_{0} \in D\left((-A)^{\epsilon}\right):\left\|u_{0}\right\|_{D\left((-A)^{\epsilon}\right)} \leq M\right\}.
$$
\begin{theorem}[see \cite{GT11}]\label{thm1}
We assume that the system \eqref{absy} is final state observable in time $T>0$. Then there exist constants $\alpha \in (0,1)$ and $ K>0$ such that
$$\left\|u_{0}\right\| \leq \frac{K}{\left|\log \|\mathbf{C} u\|_{L^{2}(0, T ; Y)}\right|^{\alpha}}\quad  \text{for all } u_{0} \in \mathcal{I}_{\epsilon, M}.$$
\end{theorem}
The main motivation of our work is to extend the above result when the generated semigroup $\left(\mathrm{e}^{t A}\right)_{t \geq 0}$ is analytic of an arbitrary angle $\psi \in \left(0,\dfrac{\pi}{2}\right]$. Note that a self-adjoint dissipative operator generates a bounded analytic semigroup of angle $\dfrac{\pi}{2}$.

\section{Analytic case}
Assume that $A$ is the generator of a bounded analytic semigroup $(\mathrm{e}^{tA})_{t\geq 0}$ of angle $\psi \in \left(0,\dfrac{\pi}{2}\right]$. Set $$\Sigma_\psi:= \{z\in \mathbb{C} \setminus \{0\}\colon |\arg z|< \psi \},$$
and let $K\geq 1$ and $\kappa \ge 0$ two constants such that
\begin{equation*}
\|\mathrm{e}^{z A}\| \leq K \mathrm{e}^{\kappa \mathrm{Re}\,z} \quad \text{ on } \overline{\Sigma}_\psi.
\end{equation*}
To prove the main theorem, we need the general logarithmic convexity result:
\begin{theorem}[Krein-Prozorovskaya (1960)]\label{thmlc} Let $u_0\in H$. The solution of \eqref{absy} satisfies
\begin{equation*}
\|u(t)\| \leq K \mathrm{e}^{\kappa(t-T w(t))} \|u_0\|^{1- w(t)} \|u(T)\|^{w(t)}, \qquad 0\leq t\leq T,
\end{equation*}
where $w$ is the harmonic function on the strip
$$\mathcal{S}_\psi =\{t + r \mathrm{e}^{\pm i \psi} \; |\; 0 \le t \le T,\; r \ge 0\},$$
which is bounded and continuous on $\mathcal{S}_\psi$, satisfying $0 = w(0) \le w(t) < w(T) = 1$ for all $0 \le t < T$.
\end{theorem}
The idea of the proof relies on the maximum principle for subharmonic functions, and the result holds for a general Banach space.

The harmonic function $w$ is given by
$$
w(z)=\frac{\operatorname{Re}\left(h^{-1}(z)\right)}{T}, \quad h(z)=f \circ g(z), \quad g(z)=T \sin ^{2}\left(\frac{\pi z}{2 T}\right)
$$
and
$$
f(z)=\frac{T \sin \psi}{\pi} \int_{0}^{\frac{z}{T}} t^{\frac{\psi}{\pi}-1}(1-t)^{-\frac{\psi}{\pi}} \mathrm{d} t.
$$
Consequently, we can bound $w$ from below.
\begin{lemma}
The harmonic function $w$ satisfies the inequality
\begin{equation*}
w(t)\ge \frac{2}{\pi} \left(\frac{\psi}{\sin\psi}\right)^{\frac{\pi}{2\psi}} \left(\frac{t}{T}\right)^{\frac{\pi}{2\psi}}, \qquad 0<t\le T.
\end{equation*}
\end{lemma}

Now, we state our main result of logarithmic stability that extends the self-adjoint case where $\psi=\frac{\pi}{2}$ (Theorem \ref{thm1}).
\begin{theorem}[see \cite{ACM'211}]\label{thm2}
Assume that the system \eqref{absy} is final state observable in $T$. Then there exist constants $\alpha \in (0,1)$ and $K_{1}>0$ such that
$$
\left\|u_{0}\right\| \leq K_{1}\left(\frac{2 \psi \Gamma\left(\frac{2 \psi}{\pi}\right)}{\pi\left|c_{\psi} \log \|\mathbf{C} u\|_{L^{2}(0, T ; Y)}\right|^\frac{2 \psi}{\pi}}\right)^{\alpha}
$$
for all $u_{0} \in \mathcal{I}_{\epsilon, M}$, where $c_{\psi}=\dfrac{c}{\pi}\left(\dfrac{\psi}{\sin \psi}\right)^{\frac{\pi}{2 \psi}}$, and $c>0$ is constant.
\end{theorem}

\subsection{Application to Ornstein–Uhlenbeck equation}
Consider the Ornstein-Uhlenbeck equation
\begin{empheq}[left =\empheqlbrace]{alignat=2}
\begin{aligned}
&\partial_t y = \Delta y + B x\cdot \nabla y, \qquad 0<t<T , && x\in \mathbb{R}^N, \\
& y\rvert_{t=0}=y_0,   && x \in \mathbb{R}^N, \nonumber
\end{aligned}
\end{empheq}
where $B$ is a real constant $N \times N$-matrix.

The operator $A:=\Delta+B x \cdot \nabla$, with its maximal domain, generates a $C_{0}$-semigroup on $L^{2}\left(\mathbb{R}^{N}\right)$. Assuming the spectral condition $
\sigma(B) \subset \{z \in \mathbb{C}: \operatorname{Re} z<0\}$ guarantees the existence of an invariant measure $\mu$ for the Ornstein-Uhlenbeck semigroup, which is analytic of angle $\psi<\dfrac{\pi}{2}$ in general (see \cite{Chill}).

For the observability inequality, we consider the observation operator $\mathbf{C}=\mathds{1}_\omega$ and the observation region $\omega \subset \mathbb{R}^N$ satisfying \cite{LM16}:
\begin{equation}\label{geom}
\exists \delta, r>0, \forall y \in \mathbb{R}^N, \exists y' \in \omega, \quad B\left(y', r\right) \subset \omega \text { and }\left|y-y'\right|<\delta. 
\end{equation}
An example of sets satisfying \eqref{geom} in $\mathbb{R}^2$ is given by $\omega:=\bigcup\limits_{(n,m)\in \mathbb{Z}^2} B((n,m),\frac{1}{2})$, where $\delta=\sqrt{2}$ and $r=\frac{1}{2}$.

Now, we recall the following observability result on the weighted space $L^2_\mu:=L^2\left(\mathbb{R}^{N}, \d\mu\right)$.
\begin{proposition}[Beauchard-Pravda Starov (2018)]
There exists a constant $\kappa_T=\kappa_T(\omega,T)>0$ such that
\begin{equation*}
\|y(T,\cdot)\|_{L^2_\mu}^2 \leq \kappa_T \int_0^T \|y(t,\cdot)\|_{L^2_\mu(\omega)}^2 \mathrm{d} t.
\end{equation*}
\end{proposition}
Applying our abstract Theorem \ref{thm2}, we obtain the stability result:
\begin{proposition}
There exist constants $\alpha\in (0,1)$ and $ K_{1}>0$ such that
$$
\left\|y_{0}\right\|_{L_{\mu}^{2}} \leq K_{1}\left(\frac{2 \psi \Gamma\left(\frac{2 \psi}{\pi}\right)}{\pi\left|c_{\psi} \log \|y\|_{L^{2}\left(0, T ; L^{2}_\mu(\omega)\right)}\right|^\frac{2 \psi}{\pi}}\right)^{\alpha} \quad \text{for all }\|y_0\|_{H^{2\epsilon}_\mu} \le M.
$$
\end{proposition}

\subsection{A non-analytic case}
Let $s$ be a positive real number. We consider the Ornstein-Uhlenbeck equation with fractional diffusion
\begin{empheq}[left =\empheqlbrace]{alignat=2}
\begin{aligned}\label{es}
&\partial_t u = -\mathrm{tr}^s\left(-Q\nabla^2 u\right) + B x\cdot \nabla u, \qquad 0<t<T , && x\in \mathbb{R}^N, \\
& u\rvert_{t=0}=u_0(x),   && x \in \mathbb{R}^N,
\end{aligned}
\end{empheq}
where $\mathrm{tr}^s\left(-Q\nabla^2 \cdot\right)$ is the Fourier multiplier whose symbol is $\left\langle Q \xi, \xi\right\rangle^s$.

P. Alphonse and J. Bernier recently studied this fractional model in the context of null controllability and proved the observability inequality (see \cite{AB20}) in the standard space $L^{2}\left(\mathbb{R}^{N}\right)$:
\begin{theorem}
We assume that $s>\frac{1}2$ and that $\omega \subset \mathbb{R}^N$ satisfy \eqref{geom}. There exists a positive constant $\kappa_T$ such that for all $u_0 \in L^{2}\left(\mathbb{R}^{N}\right)$, we have
\begin{equation*}
\|u(T,\cdot)\|_{L^{2}\left(\mathbb{R}^{N}\right)}^2 \leq \kappa_T \int_0^T \|u(t,\cdot)\|_{L^{2}(\omega)}^2\,\d t,
\end{equation*}
where $u$ is the solution associated with \eqref{es}.
\end{theorem}

We have proven the following logarithmic convexity:
\begin{proposition}\label{proplc}
Let $s>0$ and $T>0$. There exist constants $c=c(s,T)\in (0,1]$ and $K:=K(c,T)$ such that for all $u_0\in L^2\left(\mathbb R^N\right)$ the following inequality holds
\begin{equation*}
\|u(t)\|_{L^{2}\left(\mathbb{R}^{N}\right)} \le K\|u_0\|_{L^{2}\left(\mathbb{R}^{N}\right)}^{1-c\frac{t}{T}} \|u(T) \|_{L^{2}\left(\mathbb{R}^{N}\right)}^{c\frac{t}{T}}, \quad 0\le t\le T,
\end{equation*}
where $u$ is the solution of the fractional equation \eqref{es}.
\end{proposition}
The proof draws on a representation formula of the semigroup via the Fourier transform and the use of H\"older inequality.

Set $\mathcal{A}:=-\operatorname{tr}^s\left(-Q \nabla^2 \cdot\right)+B x \cdot \nabla \cdot$. Thus, we obtain a stability estimate for the class of initial data
$$\mathcal{I}_M=\left\{u_0 \in L^2\left(\mathbb{R}^N\right), \mathcal{A} u_0 \in L^2\left(\mathbb{R}^N\right):\left\|u_0\right\|_{D(\mathcal{A})} \leq M\right\}.$$
\begin{theorem}\label{thmlogstab1}
We assume that $s>\frac{1}{2}$ and that $\omega \subset \mathbb{R}^N$ satisfy \eqref{geom}. There exist positive constants $C$ and $C_1$ depending on $(s, T, \omega, M)$ such that
\begin{equation*}
\|u_0\|_{L^{2}\left(\mathbb{R}^N\right)} \leq \frac{C}{\left|\log \left(C_1\|u\|_{H^1\left(0, T ; L^{2}(\omega)\right)}\right)\right|} \quad \text{for all } u_0 \in \mathcal{I}_M,
\end{equation*}
where $u$ solves the equation \eqref{es}.
\end{theorem}

\section{Time fractional equations}
Let $0 < \alpha \le 1$ and $T>0$. We consider
\begin{empheq}[left = \empheqlbrace]{alignat=2}
\begin{aligned}
&\partial_{t}^\alpha u(t) = A u(t), && \qquad  t\in (0, T), \\
& u(0)=u_0,\\
\end{aligned} \label{eq1}
\end{empheq}
where $A: D(A) \subset H \rightarrow H$ is a densely defined operator such that  
\begin{itemize}
    \item[(i)] $A$ is self-adjoint,
    \item[(ii)] there exists $\kappa \ge 0$ such that 
$\langle A u, u\rangle \leq \kappa \|u\|^2$ for all $u \in D(A)$,
    \item[(iii)] $A$ has compact resolvent.
\end{itemize}

The Caputo derivative $\partial_{t}^\alpha g$ is defined by
\begin{equation*}
    \partial_t^{\alpha} g(t) = \begin{cases}\displaystyle \frac{1}{\Gamma(1-\alpha)} \int^t_0 (t-s)^{-\alpha} \frac{\d }{\d s}g(s) \, \d s, & 0<\alpha<1,\\
    \dfrac{\d }{\d t}g(t), & \alpha=1.
    \end{cases}
\end{equation*}
\textbf{Backward problem:} Given $u(T)$, can we recover $u(t_0)$, $0\le t_0 <T$?

We have proven the logarithmic convexity for the fractional case.
\begin{theorem}[see \cite{CMY23}] \label{thm3}
Let $0 < \alpha \le 1$. Let $u$ be the solution to \eqref{eq1}. Then there exists a constant $M\ge 1$ such that
\begin{equation*}
\|u(t)\| \le M \|u_0\|^{1-\frac{t}{T}} \|u(T)\|^{\frac{t}{T}}, \qquad 
0\le t \le T.
\end{equation*}
Moreover, if $\kappa=0$, we can choose $M=1$.
\end{theorem}
The proof's key idea relies on the solution's spectral representation and the complete monotonicity of the Mittag-Leffler functions.

For the non-symmetric case, we can consider the following operator $$L u\, (x):= \mathrm{div}(\mathcal{A}(x)\nabla u(x)) + \mathcal{B}(x)\cdot \nabla u(x) + p(x) u(x),$$ with symmetric and uniformly elliptic principal part. However, we need to introduce an assumption on the drift term that reads as follows
\begin{equation*}\label{H}
\textbf{(H)} \quad \text{There exists a function } b\in W^{2,\infty}(\Omega) \text{ such that } \mathcal{B}=\mathcal{A}\nabla b.    
\end{equation*}
The theoretical and numerical results for the non-symmetric operator $L$ can be found in \cite{CMY'25}.

\section{Open problems}
\noindent$1.$ Can we prove similar results for coupled systems? e.g., for $0<\alpha_1<\alpha_2 <1,$ and 
\begin{empheq}[left = \empheqlbrace]{alignat=2}
\begin{aligned}
&\partial_{t}^{\alpha_1} u_1 = \Delta u_1 + u_1 + u_2, \\
&\partial_{t}^{\alpha_2} u_2 = \Delta u_2 +  u_1 + u_2 . \nonumber
\end{aligned}
\end{empheq}

\noindent$2.$ Can we prove the logarithmic convexity for the non-symmetric case without Assumption \textbf{(H)}?

\noindent$3.$ Can we prove the backward uniqueness for the fractional Cauchy problem \eqref{eq1} governed by a general analytic semigroup?

\section{Acknowledgment}
The author would like to thank Marianna Chatzakou and the organizers of the Ghent Methusalem Analysis \& PDE Seminar, for which most of this material was prepared.

\newpage



\begin{thebibliography} {99}
\bibitem{ACM'211}
Ait Ben Hassi, E. M., Chorfi, S. E., and Maniar, L.: Inverse problems for general parabolic systems and application to Ornstein-Uhlenbeck equation, \emph{Discrete Contin. Dyn. Syst. - S}, \textbf{17}, 1966--1980 (2024).
\bibitem{AB20}
Alphonse, P. and Bernier, J.: Smoothing properties of fractional Ornstein–Uhlenbeck semigroups and null-controllability, {\it Bull. Sci. Math.}, \textbf{165}, 102914 (2020).
\bibitem{Chill}
Chill, R., Fašangová, E., Metafune, G. and Pallara, D.: The sector of analyticity of the Ornstein–Uhlenbeck semigroup on $L^p$ spaces with respect to invariant measure, {\it J. Lond. Math. Soc.}, \textbf{3}, 703--722 (2005).
\bibitem{CM'24}
Chorfi, S. E. and Maniar, L.: Stability Estimates for Initial Data in General Ornstein–Uhlenbeck Equations, In ``Control Theory and Inverse Problems", Birkhäuser Cham, 2024.
\bibitem{CMY'25}
Chorfi, S. E., Maniar, L. and Yamamoto, M.: Logarithmic convexity of non-symmetric time-fractional diffusion equations, {\it Math. Meth. Appl. Sci.}, \textbf{48}, 2011--2021 (2025).
\bibitem{CMY23}
Chorfi, S. E., Maniar, L. and Yamamoto, M.: The backward problem for time-fractional evolution equations, {\it Appl. Anal.}, \textbf{103}, 2194--2212 (2023).
\bibitem{GT11}
García, G. C. and Takahashi, T.: Inverse problem and null-controllability for parabolic systems, {\it J. Inverse Ill-Posed Probl.}, \textbf{19}, 379--405 (2011).
\bibitem{LM16}
Le Rousseau, J. and Moyano, I.: Null-controllability of the Kolmogorov equation in the whole phase space, {\it J. Differ. Equ.}, \textbf{260}, 3193--3233 (2016).
\bibitem{Pa75}
Payne, L. E.: Improperly Posed Problems in Partial Differential Equations, Philadelphia: SIAM (1975).
\end{thebibliography}
\end{document}